\documentclass[11pt]{amsart}
\usepackage{amscd, amssymb}
\usepackage{tikz}
\usetikzlibrary{arrows,positioning} 

\usepackage{comment}

\usepackage{amsthm}
\usepackage{mathrsfs}
\usepackage{bbm}
\usepackage{ stmaryrd }

\usepackage{multirow}
\usepackage{makecell}

\usepackage{hyperref}
\hypersetup{
	colorlinks=true,
	linkcolor=blue,
	filecolor=magenta,  
	urlcolor=cyan,
}

\def\End{\operatorname{End}}

\def\ker{\operatorname{ker}}

\def\im{\operatorname{Im}}

\def\dim{\operatorname{dim}}

\def\Ind{\operatorname{Ind}}

\def\Hom{\operatorname{Hom}}

\def\C{\mathbb{C}}

\def\N{\mathbb{N}}

\def\G{\mathbb{G}}

\def\A{\mathbb{A}}

\def\F{\mathbb{F}}

\def\A{\mathbb{A}}

\def\OO{\mathcal{O}}
\def\KK{\mathcal{K}}

\def\JJ{\mathcal{J}}

\def\EE{\mathcal{E}}

\def\DD{\mathcal{D}}

\def\PP{\mathcal{P}}

\def\Cs{\mathcal{C}}
\def\II{\mathcal{I}}

\def\As{\mathscr{A}}
\def\Bs{\mathscr{B}}
\def\Cs{\mathscr{C}}

\def\Ds{\mathscr{D}}

\def\m{\mathfrak{m}}

\def\ol{\overline}

\def\Rep{\operatorname{Rep}}

\def\Ver{\operatorname{Ver}}

\def\Spec{\operatorname{Spec}}

\def\Vec{\operatorname{Vec}}

\def\Ext{\text{Ext}}

\def\ul{\underline}

\def\Res{\operatorname{Res}}

\def\sub{\subseteq}

\def\sVec{\text{sVec}}

\def\Sp{\operatorname{Sp}}

\def\Proj{\operatorname{Proj}}

\def\xto{\xrightarrow}

\newtheorem{thm}{Theorem}[section]
\newtheorem{cor}[thm]{Corollary}
\newtheorem{lemma}[thm]{Lemma}
\newtheorem{prop}[thm]{Proposition}
\newtheorem{conj}[thm]{Conjecture}

\theoremstyle{definition}
\newtheorem{definition}[thm]{Definition}
\theoremstyle{remark}
\newtheorem{remark}[thm]{Remark}
\newtheorem{example}[thm]{Example}
\newtheorem{question}[thm]{Question}

\numberwithin{equation}{section}

\usepackage{relsize}
\usepackage{fancyhdr}
\usepackage[all,cmtip]{xy}
\begin{document}

    \title{Semi-exact functors and the conjecture of Coulembier and Flake}

    \author{Alexander Sherman}
    
    \begin{abstract}
        We introduce the notion of semi-exact functors between abelian categories, which are characterised by either preserving or breaking exactness of sequences.  We prove semi-exact functors are always exact with respect to a canonically determined exact structure. For a symmetric tensor category $\mathscr{C}$, we define a support theory by considering all $\Bbbk$-linear symmetric monoidal, semi-exact functors $\mathscr{C}\to\operatorname{Ver}_{p^\infty}$ up to equivalence, and conjecture that it agrees with the Balmer spectrum of the stable category when $\mathscr{C}$ is finite.  We prove a weak form of this conjecture and other results for Tannakian categories and in characteristic 0.  Critical to the positive characteristic case is a functor constructed by Coulembier and Flake who were interested in higher Frobenius functors.  We show that their construction, which relies on the theory of abelian envelopes, gives a natural source for semi-exact functors, and at least in the Tannakian setting gives all $\Bbbk$-linear, symmetric monoidal, semi-exact functors.
    \end{abstract}

\maketitle

\section{Introduction: semi-exact functors}

Throughout, $\Bbbk$ denotes an algebraically closed field.  We begin with a definition.

\begin{definition}
    Let $\As,\Bs$ be abelian categories.  We say an additive functor $\varphi:\As\to\Bs$ is \textbf{semi-exact} if for any short exact sequence
    \[
    0\to A\to B\to C\to 0,
    \]
    the following are equivalent:
    \begin{enumerate}
        \item the sequence
    \[
    0\to \varphi(A)\to\varphi(B)\to\varphi(C)\to0
    \]
    is exact;
    \item the sequence
    \[
    0\to \varphi(A)\to\varphi(B)
    \]
    is exact;
    \item the sequence 
    \[
    \varphi(B)\to\varphi(C)\to0
    \]
    is exact.
    \end{enumerate}
\end{definition}

It is clear that any exact functor is semi-exact; however, there are examples of semi-exact functors which are not exact.  For example, if $\As=A$-mod, where $A=\Bbbk[x]/x^2$, and $\Bs=\Vec$, then the functor $\ker(x)/\operatorname{im}(x)$ is semi-exact.

The notion of semi-exact functor does not appear in the literature to the best of this author's knowledge. The purpose of this paper is to make a study of such functors.  We find that in the setting of symmetric tensor categories, semi-exact functors have remarkable properties, and a rich connection to support theory.

\subsection{Conjectures and theorems}  By a symmetric tensor category $\Cs$ over $\Bbbk$ we mean a symmetric monoidal, rigid, $\Bbbk$-linear, locally finite abelian category such that $\End(\mathbf{1})=\Bbbk$.  If $\Cs$ has finitely many simples and enough projectives, we additionally say that $\Cs$ is finite.

If $\Cs$ is a finite symmetric tensor category, then we may form its stable module category $\underline{\Cs}$, which is a tensor triangulated category, and therefore has a Balmer spectrum $\Sp\underline{\Cs}$ consisting of all thick prime tensor ideals in $\underline{\Cs}$.  We write $\operatorname{cl}\Sp\underline{\Cs}$ for the closed points of the topological space $\Sp\underline{\Cs}$.

In the following conjecture, let $p$ denote the characteristic of our field $\Bbbk$, which is either a prime number or $0$.  If $p=0$, we set $\Ver_{p^\infty}=\sVec$, and for a prime $p>0$ we let $\Ver_{p^\infty}$ denote the union $\bigcup\limits_{n}\Ver_{p^n}$ (\cite{BEO}).  We use the abbreviation LSM functor for a $\Bbbk$-linear, symmetric monoidal functor (not necessarily exact).

\begin{conj}\label{conj 1}
    Let $\Cs$ be a finite symmetric tensor category.  Then there is a bijection:

\[
\left\{
\begin{array}{c}
    \text{non-exact, semi-exact LSM functors}  \\
      \varphi:\Cs\to\Ver_{p^\infty}
\end{array}
\right\}_{\bigg/\sim} \longrightarrow \operatorname{cl}\operatorname{Sp}\ul{\Cs},
\]
where on the left hand side we have $\varphi\sim\psi$ if the two functors $\varphi,\psi:\Cs\to\Ver_{p^\infty}$ are naturally isomorphic as monoidal functors.

\end{conj}

Let us explain where the map in Conjecture \ref{conj 1} comes from.  Given a non-exact, semi-exact LSM functor $\varphi:\Cs\to\Ver_{p^\infty}$, its object kernel will define a point $x\in\Sp\ul{\Cs}$ such that $\varphi(M)\neq 0$ if and only if $x$ lies in the support of $M$.  If we are able to identify $\Sp\ul{\Cs}$ with $\Proj \Ext^\bullet(\mathbf{1},\mathbf{1})$ (which happens for Tannakian categories \cite{BIKP} and conjecturally always holds \cite{NVY}), then we show that $x$ must be a closed point.  This gives a map in one direction.  However it may be surprising that this correspondence could be a bijection--we return to this shortly.  In the meanwhile, here is another (perhaps more outlandish) conjecture.  We will say that an LSM functor $\varphi:\Cs\to\Ds$ is surjective if every object in $\Ds$ is a subquotient of an object in the image of $\varphi$.

\begin{conj}\label{conj 2}
     Let $\Cs,\Ds$ be symmetric tensor categories, and additionally assume that $\Cs$ is finite.  Let $\varphi:\Cs\to\Ds$ be a surjective, semi-exact, LSM functor.  Then $\Ds$ is also finite.
\end{conj}

Note that Conjecture \ref{conj 2} holds for exact functors \cite[Prop.~1.8.19]{EGNO}, but is false if we drop the semi-exactness property; indeed, semisimplification functors are a large source of counterexamples.  Here is some immediate evidence for Conjectures \ref{conj 1} and \ref{conj 2}.

\begin{thm}\label{thm intro}
    Conjectures \ref{conj 1} and \ref{conj 2} hold if $\Bbbk$ is characteristic 0 or $\Cs$ is Tannakian, i.e.~$\Cs\simeq\Rep G$ for a finite group scheme $G$.
\end{thm}

\subsection{Connection with abelian envelopes}\label{intro ab envelopes} What is perhaps interesting about the proof of Theorem \ref{thm intro} is it can be done with a uniform argument that suggests how the proofs of Conjectures \ref{conj 1} and \ref{conj 2} would go in general.  The key observation is that in each case, given a closed point $\PP\in\operatorname{cl}\Sp\ul{\Cs}$, we can construct a nonzero additive tensor ideal $\II_\PP$ with finitely many indecomposables that are all supported at the point $\PP$.  From this, one can turn the crank of the theory of abelian envelopes (see \cite[Thm.~2.4.1]{CEO}, \cite{St}) to obtain a functor
\[
\varphi_{\PP}:\Cs\to (\ul{\Cs}/\PP^{max})^{ab},
\]
where $\PP^{max}$ is a certain tensor ideal in $\ul{\Cs}$, $(\ul{\Cs}/\PP^{max})^{ab}$ is the abelian envelope of $\ul{\Cs}/\PP^{max}$, and is a finite symmetric tensor category.  We prove that the functor $\varphi_\PP$ is semi-exact, and moreover in good situations any semi-exact functor on $\Cs$ must factor through one of the form $\varphi_\PP$  for some $\PP\in\operatorname{cl}\Sp\ul{\Cs}$.  Using the recent work \cite{CF2} which proves the Benson-Etingof-Ostrik conjecture for finite tensor categories, along with uniqueness of fibre functors to $\Ver_{p^\infty}$, Theorem \ref{thm intro} follows.

\subsection{The conjecture of Coulembier-Flake}  While Section \ref{intro ab envelopes} gives a path to proof, the very difficult part is to construct an additive tensor ideal $\II_\PP$ with the prescribed properties.  In characteristic 0 this is not so difficult, but in characteristic $p>0$ we require (a generalisation of) the `Saturation Conjecture' of Coulembier-Flake \cite[Conj.~5.4.2(1)]{CF}.  This conjecture provides such additive tensor ideals in the case of quasi-elementary abelian groups, meaning ones isomorphic to $C_p^r\times\G_{a(s)}$, where $\G_{a(s)}$ denotes the $s$th Frobenius kernel of $\G_a$, and $C_p$ is the cyclic group of order $p>0$.  Given such a quasi-elementary abelian group $E$ of rank $n=r+s$, along with a generic point $\PP\in\operatorname{cl}\Sp\ul{\Rep E}$, \cite{CF} provides an additive tensor ideal $\II_\PP$ with the desired properties, and a surjective, semi-exact LSM functor
\[
\varphi_\PP:\Rep E\to\Ver_{p^n}.
\]
The Saturation Conjecture of Coulembier-Flake was recently proven in \cite{CF2}.  We will provide a different more geometric proof which was discovered independently.

\subsection{Upshot in characteristic 0 setting}  The upshot of Theorem \ref{thm intro} in characteristic 0 is as follows. In this case, Deligne's theorem \cite{Del} implies that $\Cs\simeq\Rep(G,\epsilon)$ for a supergroup $G$ of the form $H\ltimes\Pi V$, where $H$ is a finite group and $V$ is an $H$-module. From an element $v\in V$ we obtain a Duflo-Serganova functor (\cite{DS}, \cite{GHSS})
\[
DS_v:\Cs\to\sVec,
\]
which is a semi-exact (by Hinich's lemma, for the experts in superalgebras) LSM functor.  Theorem \ref{thm intro} implies that any semi-exact functor on $\Cs$ must be isomorphic to $DS_v$ for some $v\in V$.

\subsection{Upshot in characteristic $p>0$}  Let $G$ be a finite group scheme in characteristic $p>0$, and write $\Cs=\Rep G$.  Then every point $x\in\operatorname{cl}\Sp\ul{\Cs}$ arises as a generic closed point in $\operatorname{cl}\Sp\ul{\Rep E}$, where $E\leq G$ is a quasi-elementary abelian subgroup of $G$ \cite{FP}.  Using the construction of Coulembier-Flake, we obtain, by restriction to $E$, a semi-exact LSM functor
\[
\varphi_x:\Rep G\to\Ver_{p^n}\sub\Ver_{p^\infty}.
\]
The upshot of Theorem \ref{thm intro} is that every (non-exact) semi-exact LSM functor on $\Rep G$ arises in this way.

\subsection{Connection to support theory}  The bijection in Conjecture \ref{conj 1} point toward a purely abelian approach to support theory which may be useful even beyond the finite setting.  

For a field extension $\mathbb{L}/\Bbbk$, define $(\Ver_{p^\infty})_{\mathbb{L}}$ to be the extension of scalars of $\Ver_{p^\infty}$, which is once again a symmetric tensor category \cite[Lem.~2.2]{EO}.  If $\mathbb{L},\mathbb{L}'$ are field extensions of $\Bbbk$, and
\[
\varphi_1:\Cs\to(\Ver_{p^\infty})_{\mathbb{L}} \ \text{ and } \ \varphi_2:\Cs\to(\Ver_{p^\infty})_{\mathbb{L}'}
\]
are LSM functors, then we write $\varphi_1\sim\varphi_2$ if there exists a field extension $\mathbb{L}''$ of both $\mathbb{L}$ and $\mathbb{L}'$ such that $\varphi_1$ and $\varphi_2$ are coequalised by the natural functors $(\Ver_{p^\infty})_{\mathbb{L}}\to(\Ver_{p^\infty})_{\mathbb{L}''}$ and $(\Ver_{p^\infty})_{\mathbb{L}'}\to(\Ver_{p^\infty})_{\mathbb{L}''}.$  This defines an equivalence relation on the set of LSM functors $\Cs\to(\Ver_{p^\infty})_{\mathbb{L}}$ ranging over all field extensions of $\mathbb{L}/\Bbbk$.  Now define
\[
SE(\Cs)=\{\varphi:\Cs\to(\Ver_{p^\infty})_{\mathbb{L}}:\varphi\text{ is a non-exact, semi-exact LSM functor, } \mathbb{L}/\Bbbk\text{ a field ext}\}/\sim.
\]
We equip $SE(\Cs)$ with a topology by declaring a basis for the closed sets to be those of the form $I(X)=\{\varphi:\varphi(X)\neq0\}$, where $X$ runs over all objects of $\Cs$.  Then for $X\in\Cs$, we define its semi-exact support to be
\[
\operatorname{supp}_{SE}(X):=I(X).
\]
By definition $\operatorname{supp}_{SE}(X)$ is closed, and in this way $(SE(\Cs),\operatorname{supp}_{SE})$ defines a support theory  on $\Cs$ which can easily be shown to satisfy the tensor product property:
\[
\operatorname{supp}_{SE}(X\otimes Y)=\operatorname{supp}_{SE}(X)\cap\operatorname{supp}_{SE}(Y).
\]
One advantage of $\operatorname{SE}(\Cs)$ is that it does not require $\Cs$ to be quasi-Frobenius. However, if $\Cs$ is quasi-Frobenius, then $(SE(\Cs),\operatorname{supp}_{SE})$ defines a support theory on $\ul{\Cs}$.

\begin{conj}\label{conj 3}
    If $\Cs$ is a finite symmetric tensor category, then $\operatorname{SE}(\Cs)$ is homeomorphic to $\Sp\ul{\Cs}$, and $(SE(\Cs),\operatorname{supp}_{SE})$ is the universal support theory for $\ul{\Cs}$.
\end{conj}

It is possible that Conjecture \ref{conj 3} holds for any quasi-Frobenius category, but we have not computed any examples yet.

\subsection{Statement on AI usage and story} This project began when the author discovered a proof of the Coulembier-Flake Saturation Conjecture using ChatGPT Astra 6.  While writing up this result, the author became aware that Coulembier and Flake had recently also found a proof, which can now be found in \cite{CF2} (see also the proof given by ChatGPT in \cite{chatgpt}). However the author was originally motivated by the semi-exactness property in connection with the paper \cite{HKPSYZ}, which led to this paper.

The author used ChatGPT Astra 6 to help with the proof of Theorem \ref{thm intro}, whose statement evolved over the last two months of work.  The author also used AI to find references, understand related works in the literature, and to check the final version.  Anyone interested in a more detailed AI usage statement is welcome to write the author.

\textbf{This paper has been written entirely by the author, and the author takes full responsibility for its accuracy.}

\subsection{Acknowledgements} The author would like to foremost thank Kevin Coulembier and Johannes Flake for their comments on an initial draft of this preprint, and for helpful discussions thereafter.  The author gives a hearty thanks to Bregje Pauwels for many insightful discussions and for coining the term semi-exact.  The author further thanks Jonathan Gruber for interesting conversations about this project. Finally, the author thanks OpenAI for providing a temporary ChatGPT for Academic Researchers account.

\section{Semi-exact functors, continued}

We return to the study of semi-exact functors $\varphi:\As\to\Bs$, where $\As,\Bs$ are now arbitrary abelian categories.  Let us begin with an example.

\begin{example}\label{example 1}
    Consider $A=\Bbbk[x]/(x^n)$, and let $\As=A$-mod, $\Bs=\Vec$.  For $1\leq i\leq n$, define $\varphi_i:\As\to\Bs$ by 
    \[
\varphi_i=\frac{\ker(x)\cap\operatorname{im}(x^{i-1})}{{\ker(x)\cap\operatorname{im}(x^i)}}.
    \]
    Then $\bigoplus\limits_{i=1}^{n-1}\varphi_i$ is a semi-exact functor \cite[Lem.~3.3.3]{HKPSYZ}.  

    In particular, in characteristic $p>0$ the semisimplification functor $\Rep C_p\to\Ver_p$, which is given by
    \[
    \bigoplus\limits_{i=1}^{p-1}\varphi_i\boxtimes L_i,
    \]
    is a semi-exact functor.  (For the definition and properties of $\Ver_p$, we refer to \cite{EO}.)
\end{example}

\begin{remark}
    Example \ref{example 1} shows that a semi-exact functor need not be middle-exact. 
\end{remark}

Here is a simple but useful lemma, whose proof we omit.

\begin{lemma}\label{lemma easy}
    Suppose that $\varphi:\As\to\Bs$ is semi-exact, and $0\to X\to Y\to Z\to 0$ is a short exact sequence in $\As$.  
    \begin{enumerate}
        \item If $\varphi(X)=0$, then $\varphi(Y)\cong \varphi(Z)$.  
        \item If $\varphi(Z)=0$, then $\varphi(X)\cong\varphi(Y)$.
    \end{enumerate}
\end{lemma}

\subsection{Exact structure from a semi-exact functor}Let $\varphi:\As\to\Bs$ be a semi-exact functor     between abelian categories.  Consider the class $\EE_{\varphi}$ of short exact sequences $0\to A\to B\to C\to 0$ in $\As$ such that $0\to\varphi(A)\to\varphi(B)\to\varphi(C)\to0$ is exact.  
\begin{prop}
  The set $(\As,\EE_{\varphi})$ is an exact category, and $\varphi:\As\to\Bs$ is exact with respect to this exact structure.
\end{prop}

\begin{proof}
    We verify the axioms of an exact category.  The only nontrivial axioms to check are that the pushout of admissible monomorphisms and the pullback of admissible epimorphisms remain admissible.  We will prove it for admissible monomorphisms, with the case of admissible epimorphisms being similar.  Let $0\to A\to B$ be a monomorphism in $\As$ such that $0\to\varphi(A)\to\varphi(B)$ is exact.  Consider a morphism $A\to C$ giving the pushforward $D$ which fits into a short exact sequence
    \[
    0\to A\to B\oplus C\to D\to 0.
    \]  
    Because $A\to B$ is admissible, $\varphi(A)\to\varphi(B\oplus C)=\varphi(B)\oplus\varphi(C)$ will be a monomorphism, and thus by definition of a semi-exact functor, the following sequence is exact:
    \[
    0\to \varphi(A)\to\varphi(B)\oplus\varphi(C)\to\varphi(D)\to 0.
    \]
    It follows that $\varphi(C)\to \varphi(D)$ is the pushout of $\varphi(A)\to\varphi(B)$ along $\varphi(A)\to\varphi(C)$.  Since the pushout of a monomorphism is a monomorphism in an abelian category, $\varphi(C)\to\varphi(D)$ is a monomorphism, implying that $C\to D$ is an admissible monomorphism, as desired.

    It is now obvious that $\varphi$ is exact on $(\As,\EE_{\varphi})$, so this completes the proof.
\end{proof}

\subsubsection{A further exactness property}
We say that an abelian category $\As$ is \textbf{quasi-Frobenius} if it has enough injectives, enough projectives, and injectives and projectives agree.  If $\As$ is quasi-Frobenius, we can consider the stable category $\ul{\As}$, which is the quotient of $\As$ by the ideal of morphisms factoring through a projective (equivalently, an injective) object.  The category $\ul{\As}$ has the structure of a triangulated category, where the triangulated shift is given by
\[
\Sigma(X)=X[1]\cong I/X,
\]
where $I$ is the injective hull of $X$.  Moreover, every exact triangle in $\ul{\As}$ is of the form
\[
X\xto{\alpha} Y\xto{\beta} Z\xto{\delta} \Sigma( X),
\]
where $0\to X\xto{\alpha} Y\xto{\beta} Z\to 0$ is a short exact sequence in $\As$. We do not describe the boundary map $\delta$ here.  

\begin{lemma}\label{lemma semi exact q frob}
    Suppose that $\As$ is quasi-Frobenius, $\varphi:\As\to\Bs$ is semi-exact, and $\varphi(P)=0$ for all projectives $P\in\As$.  Let
    \[
    0\to X\to Y\to\ Z\to 0
    \]
    be a short exact sequence such that $\varphi(Y)=0$.  Then we have a natural isomorphism
    \[
   \varphi(Z)\cong\varphi(\Sigma(X)).
    \]
\end{lemma}

\begin{proof}
    We have an exact triangle $Y\to Z\to \Sigma(X)\to \Sigma(Y)$ in $\ul{\As}$.  Thus there exists a short exact sequence $0\to Y\to Z\oplus I\to \Sigma(X)\to 0$ in $\As$, for some injective $I$ in $\As$.  Since $\varphi(Y)=0$, the semi-exactness property implies that $\varphi(Z)\cong \varphi(Z\oplus I)\cong \varphi(\Sigma(X))$ as desired.  
\end{proof}

\subsection{Preliminaries on symmetric monoidal categories}  A symmetric pseudo-tensor category $\Cs$ will mean a Karoubian, $\Bbbk$-linear, rigid symmetric monoidal category.  We recall that a symmetric tensor category is a symmetric pseudo-tensor category which is additionally abelian and locally finite.

Given a symmetric pseudo-tensor $\Cs$, a \textbf{tensor ideal} $\KK$ of $\Cs$ is an ideal of $\Cs$ (meaning a choice of subgroups $\KK(X,Y)\leq\Cs(X,Y)$ satisfying the usual closure properties) such that $\operatorname{id}_X\otimes f\in\KK$ for all $f\in\KK$.  

An \textbf{additive tensor ideal} $\II$ of $\Cs$ (also known as a thick tensor ideal in some places) is a collection of objects of $\Cs$, closed under isomorphisms and direct sum(mands), such that if $X\in\II$ and $Y\in\Cs$, then $X\otimes Y\in\II$.

Given a tensor ideal $\KK$ of $\Cs$, we write $\operatorname{Ob}(\KK)$ for the collection of objects $X\in\Cs$ such that $\operatorname{id}_X\in\KK$, which will be an additive tensor ideal. 

Given an LSM functor $\varphi:\Cs\to\Ds$ of symmetric pseudo-tensor categories, we will write $\ker\varphi$ for tensor ideal of morphisms $f\in\Cs$ such that $\varphi(f)=0$, and set $\KK_{\varphi}=\operatorname{Ob}(\ker \varphi)$.  

\subsection{Semi-exact LSM functors}

If $\Cs$ is a symmetric tensor category over $\Bbbk$ which has enough projectives, then it will be quasi-Frobenius and we may define its stable module category $\ul{\Cs}$.  In this case $\ul{\Cs}$ will be a symmetric monoidal, rigid, $\Bbbk$-linear triangulated category (otherwise known as a tensor triangulated category).

The following lemma is well-known \cite[Lem.~2.2.2]{CF}.
\begin{lemma}\label{lemma exactness projs}
    If $\varphi:\Cs\to\Ds$ is an LSM functor and $\Cs$ is quasi-Frobenius, then the following are equivalent:
    \begin{enumerate}
        \item $\varphi$ is exact;
        \item $\varphi(P)\neq0$ for all non-zero projectives $P$ in $\Cs$;
        \item $\varphi(P)\neq0$ for some non-zero projective $P$ in $\Cs$.
    \end{enumerate}
\end{lemma}

Let $\Cs$ be quasi-Frobenius, and let $\varphi:\Cs\to\Ds$ be a non-exact, semi-exact LSM functor.  By Lemma \ref{lemma exactness projs}, $\varphi(P)=0$ for all projectives $P$ in $\Cs$, and thus $\varphi$ factors through a functor $\ul{\Cs}\to\Ds$, which by abuse of notation we also denote $\varphi$.  Write 
\[
L:=\varphi(\Sigma(\mathbf{1})).
\]
Since $\Sigma(\mathbf{1})$ is invertible in $\ul{\Cs}$, $L$ is invertible in $\Ds$, and thus is simple.  Moreover, we have 
\[
\varphi(\Sigma^i(\mathbf{1}))\cong\varphi(\Sigma(\mathbf{1}))^{\otimes i}\cong L^{\otimes i}.
\]
\begin{lemma}
    If $0\to X\to Y\to Z\to 0$ is a short exact sequence in $\Cs$ such that $\varphi(Y)=0$, then we have a natural isomorphism
    \[
    \varphi(Z)\cong\varphi(X)\otimes L.
    \]
\end{lemma}
\begin{proof}
    This follows from Lemma \ref{lemma semi exact q frob} and the fact that $\varphi$ is monoidal.
\end{proof}

\subsubsection{The Balmer spectrum} We refer to \cite{Bal} for more details of the following construction.  We continue to assume that $\Cs$ is a quasi-Frobenius symmetric tensor category.  Recall that a \textbf{prime thick tensor ideal} in $\ul{\Cs}$ is a proper triangulated full subcategory $\PP\sub\ul{\Cs}$ such that
\begin{enumerate}
    \item $X\otimes Y\in\PP$ if $X\in \PP$ and $Y\in\ul{\Cs}$, and
    \item if $X\otimes Y\in\PP$, then either $X\in\PP$ or $Y\in\PP$.
\end{enumerate}
We define the Balmer spectrum $\operatorname{Sp}\ul{\Cs}$ to be the topological space whose points are thick prime tensor ideals in $\ul{\Cs}$, with a basis of closed sets given by $I(X)=\{\PP\in\Sp\ul{\Cs}:X\notin\PP\}$ where $X$ runs over objects of $\Cs$. We write $\operatorname{cl}\Sp\ul{\Cs}$ for the closed points of $\Sp\ul{\Cs}$.  Recall that given an object $X\in\ul{\Cs}$, its \textbf{support} is defined to be $I(X)$.

\begin{lemma}\label{lemma semi exact gives prime ideal}
    Let $\varphi:\Cs\to\Ds$ be a non-exact, semi-exact, LSM functor.  If $\Cs$ is quasi-Frobenius, then $\KK_{\varphi}=\operatorname{Ob}(\ker\varphi)$ descends to a prime thick tensor ideal of $\ul{\Cs}$. 
\end{lemma}

\begin{proof}
   If $\varphi$ is not exact, then $\KK_{\varphi}$ must contain all projective objects, so $\KK_{\varphi}$ descends to an additive tensor ideal of $\ul{\Cs}$.  To show $\KK_{\varphi}$ is a prime thick tensor ideal, we first observe that it clearly satisfies properties (1) and (2), so it remains to show that it is a triangulated full subcategory of $\ul{\Cs}$.  For this, it suffices to prove that if $X\to Y\to Z\to X[1]$ is a triangle and $X,Y\in\KK_{\varphi}$, then $Z\in\KK_{\varphi}$.  However all triangles in $\ul{\Cs}$ arise from short exact sequences in $\Cs$, thus we may apply e.g.~Lemma \ref{lemma easy} for the result.
\end{proof}

\begin{cor}\label{cor support}
    Let $\varphi:\Cs\to\Ds$ be a non-exact, semi-exact LSM functor.  Then there exists $\PP\in\Sp\ul{\Cs}$ such that for all $M\in\Cs$, $\varphi(M)\neq0$ if and only if $\PP\in\operatorname{supp}(M)$.
\end{cor}

Lemma \ref{lemma semi exact gives prime ideal} gives a mapping:
\[
\kappa:\left\{
\begin{array}{c}
    \text{non-exact, semi-exact LSM functors}  \\
      \varphi:\Cs\to\Ver_{p^\infty}
\end{array}
\right\}/_\sim \longrightarrow\Sp\ul{\Cs}.
\]
\begin{question}\label{question}
    Is $\kappa$ a bijection onto $\operatorname{cl}\Sp\ul{\Cs}$?
\end{question}

\begin{remark}
    Note that Conjecture \ref{conj 1} is the case of Question \ref{question} when $\Cs$ is additionally assumed to be finite.  A natural first case to consider beyond the finite setting is representations of a quasi-reductive supergroup in characteristic 0.
\end{remark}

\subsection{The $\Ext$-algebra}
Let $\Cs$, $\Ds$ be symmetric tensor categories, and suppose that $\Cs$ is quasi-Frobenius.  Let $\varphi:\Cs\to\Ds$ denote a non-exact, semi-exact LSM functor, and write $L:=\varphi(\Sigma(\mathbf{1}))$.  Consider the $\N$-graded $\Bbbk$-algebra 
\[
R_{L}=\bigoplus\limits_{i\geq 0}\Hom(\mathbf{1},L^{\otimes i}),
\] 
where multiplication is defined as follows.  Given $f:\mathbf{1}\to L^{\otimes i}$, $g:\mathbf{1}\to L^{\otimes j}$, define $gf$ to be 
\[
\mathbf{1}\xto{f}L^{\otimes i}\xto{g\otimes L^{\otimes i}}L^{\otimes (i+j)}.
\]
\begin{lemma}
    $R_{L}$ is a domain.
\end{lemma}
\begin{proof}
    Indeed, since $L^{\otimes i}$ is simple for all $i$, any non-zero $f:\mathbf{1}\to L^{\otimes i}$ must be an isomorphism, so any product of non-zero homogeneous elements must be non-zero. 
\end{proof}

Now write $R_{\Cs}=\bigoplus\limits_{n\geq 0}\Ext^n(\mathbf{1},\mathbf{1})$, which is a graded-commutative $\Bbbk$-algebra \cite[Prop.~3.3]{Bal}.  We obtain an algebra morphism
\[
r_{\varphi}:R_{\Cs}\to R_L,
\]
via the maps
\[
\Ext^i(\mathbf{1},\mathbf{1})\cong \Hom_{\ul{\Cs}}(\mathbf{1},\Sigma^i(\mathbf{1}))\xto{\varphi}\Hom_{\Ds}(\mathbf{1},\varphi(\Sigma^i(\mathbf{1})))\cong\Hom_{\Ds}(\mathbf{1},L^{\otimes i}).
\]
Let $\m_{\varphi}=\{\psi\in R_{\Cs}:r_{\varphi}(\psi)=0\}$ denote the kernel of $r_\varphi$.  Then since $R_{L}$ is a domain, $\m_{\varphi}$ is a homogeneous prime ideal of $R_{\Cs}$.

Let us write $\Spec^h(R_{\Cs})$ for the homogeneous prime spectrum of $R_{\Cs}$, meaning the set of homogeneous prime ideals of $R_{\Cs}$.  For $\PP\in\Sp\ul{\Cs}$, define
\[
\rho(\PP)=\sum\limits_{i}\{f:\mathbf{1}\to\Sigma^i(\mathbf{1}):\operatorname{cone}(f)\notin\PP\}.
\]
Then $\rho(\PP)$ defines a prime ideal in $R_{\Cs}$ \cite[Thm.~5.3]{Bal}.  Hence we obtain a map
\[
\rho:\Sp\ul{\Cs}\to\Spec^h(R_{\Cs}).
\]
Observe that $\operatorname{Proj}(R_{\Cs})$ is a subspace of $\Spec^h(R_{\Cs})$.  We now introduce the following assumption:
\begin{center}
    \textbf{Assumption ($\star$)}: $\rho$ is a homeomorphism of $\Sp\ul{\Cs}$ onto $\operatorname{Proj}R_{\Cs}.$
\end{center}

Note that assumption ($\star$) is known to hold for finite symmetric tensor categories in characteristic 0 (\cite{BKN}, \cite{HV} and \cite{DK}), for representations of finite group schemes \cite{BIKP}, and for representations of unipotent finite supergroup schemes \cite{BIKP2}.  Conjecturally it is expected to hold for all finite symmetric tensor categories \cite{NVY}.

\begin{thm}\label{thm m_phi maximal}
   Continuing with our assumptions above, we have the following:
    \begin{enumerate}
        \item Either $\m_{\varphi}$ is the irrelevant ideal, or $\m_{\varphi}$ is a closed point in $\operatorname{Proj}R_{\Cs}$ .
        \item If assumption ($\star$) holds, then under this identification $\KK_{\varphi}=\m_{\varphi}$, and thus $\KK_{\varphi}$ is a closed point.
    \end{enumerate}
\end{thm}

\begin{proof}
    For part (1), suppose that $\m_{\varphi}$ is not the irrelevant ideal. Then $\dim\Hom_{\Ds}(\mathbf{1},L^{\otimes i})\leq 1$.  It follows that $\im r_{\varphi}$ has Krull dimension 1, and so $\operatorname{Proj}(\im r_{\varphi})$ is a point, and $\m_{\varphi}$ will be the pullback of this point under $r_{\varphi}$, and thus is closed.

    For part (2), we need to prove that for $f:\mathbf{1}\to\Sigma^i(\mathbf{1})$, $\operatorname{cone}(f)\notin\KK_{\varphi}$ if and only if $f\in\m_{\varphi}$. And indeed, we see that $r_{\varphi}(f)\neq 0$ if and only if $0\to \mathbf{1}\xto{\varphi(f)}L^{\otimes i}$ is exact, which by semi-exactness is equivalent to
    \[
    0\to \mathbf{1}\xto{\varphi(f)}L^{\otimes i}\to\varphi(\operatorname{cone}(f))\to0
    \]
    being exact.  However $L^{\otimes i}$ is simple, meaning $\varphi(f)$ must be an isomorphism, so this is equivalent to $\varphi(\operatorname{cone}(f))=0$.  Now Corollary \ref{cor support} implies this is equivalent to $\operatorname{cone}(f)\notin \KK_{\varphi}$, as desired. 
\end{proof}

\begin{remark}
    It is possible that $\m_{\varphi}$ is the irrelevant ideal if $\Cs$ is not finite.  For example, consider the supergroup $G=\G_m\ltimes\Pi\Bbbk$, where $\G_m$ acts nontrivially on $\Pi\Bbbk$, and set $\Cs=\Rep(G)$.  Then $R_{\Cs}=\Bbbk$, meaning every prime ideal is irrelevant.
\end{remark}

\subsection{Saturation, abelian envelopes, and semi-exact functors}\label{sec abel env}

Let $\Cs$ be a quasi-Frobenius symmetric tensor category, and let $\PP\in\operatorname{cl}\Sp\ul{\Cs}$.    We say that a non-zero additive tensor ideal $\II$ of $\ul{\Cs}$ \textbf{realises} $\PP$ if
\begin{enumerate}
    \item every nonzero object $X\in\II$ has support precisely $\PP$, i.e.~$X\notin\PP$ but $X\in\PP$ for all $\PP\in\Sp\ul{\Cs}\setminus\{\PP\}$, and
    \item $\II$ has finitely many indecomposable objects.
\end{enumerate}

Suppose that $\II$ is an additive tensor ideal in $\ul{\Cs}$ which realises $\PP$. Let $\PP^{max}$ denote the largest tensor ideal of $\ul{\Cs}$ such that $\operatorname{Ob}(\PP^{max})=\PP$.  By an abelian envelope of $\ul{\Cs}/\PP^{max}$ we mean a symmetric tensor category $\Ds$ and a faithful LSM functor $\ul{\Cs}/\PP^{max}\to\Ds$ such that if $\Ds'$ is another symmetric tensor category then any faithful LSM functor $\ul{\Cs}/\PP^{max}\to\Ds'$ factors uniquely through an exact LSM functor $\Ds\to\Ds'$.

\begin{thm}\label{thm abelian envelope semi exact}
    There exists an abelian envelope $\ul{\Cs}/\PP^{max}\to\Ds$, where $\Ds$ is finite, such that the induced functor $\varphi:\Cs\to\Ds$ given by the composition
    \[
    \Cs\to\ul{\Cs}\to\ul{\Cs}/\PP^{max}\to\Ds,
    \]
    is semi-exact.
\end{thm}

For the proof, we begin with a lemma.
\begin{lemma}\label{lemma unique cover}
     There exists an additive tensor ideal $\JJ\leq\II$ such that $\JJ+\PP$ is the unique minimal additive tensor ideal in $\ul{\Cs}$ which contains $\PP$.  In other words if $\JJ'$ is another additive tensor ideal in $\ul{\Cs}$ strictly containing $\PP$, then it must contain $\JJ+\PP$.  
\end{lemma}

\begin{proof}
    We take for $\JJ$ a minimal non-zero additive thick tensor ideal in $\II$, which exists because $\II$ has finitely many indecomposables.  Then suppose that $\JJ'\leq \ul{\Cs}$ is an additive tensor ideal which strictly contains $\PP$.  Then since $\JJ'\neq\PP$, there must exist $X\in\JJ'$ which is supported at $\PP$.  Therefore for non-zero $Y\in\JJ$ we must have $X\otimes Y\in\JJ\cap\JJ'$, meaning that $\JJ'\cap\JJ$ is a non-zero additive tensor ideal contained in $\JJ$. Since $\JJ$ is minimal, it must be equal to $\JJ$, and therefore $\JJ'$ must contain $\JJ+\PP$.
\end{proof}

\begin{proof}[Proof of Theorem \ref{thm abelian envelope semi exact}]
    Let $\JJ$ be the additive tensor ideal from Lemma \ref{lemma unique cover}.  Then $\PP$ satisfies the assumptions of \cite[Thm.~2.4.1 (2)]{CEO}, where the minimal cover is $\JJ+\PP$.  Thus $\ul{\Cs}/\PP^{max}$ admits an envelope $\psi:\ul{\Cs}/\PP^{max}\to\Ds$, and moreover (1) $\psi$ restricted to $\JJ$ is full and (2) $\psi(\JJ)$ will be the category of projectives in $\Ds$.

    It remains to prove that $\varphi:\Cs\to\Ds$ is semi-exact.  Let $0\to A\to B\to C\to 0$ be a short exact sequence in $\Cs$, and suppose that $0\to \varphi(A)\to\varphi(B)$ is exact.  Let $P\in\JJ$ be non-zero.  Tensoring with $\varphi(P)$, we obtain an exact sequence
    \[
    0\to \varphi(A)\otimes\varphi(P)\to\varphi(B)\otimes\varphi(P).
    \]
    However $\varphi(P)$ is projective, and thus injective, meaning each term in the exact sequence is injective.  Hence we obtain a retraction $r:\varphi(B\otimes P)\to\varphi(A\otimes P)$.  Because $\varphi$ is fully faithful on $\JJ$, there exists a map $r':A\otimes P\to B\otimes P$ in $\Cs$ such that $r'\circ (P\otimes \alpha)=1+n$, where $n\in\PP^{max}$.  However $\PP^{max}$ does not contain the identity on any object in $\JJ$, and thus $n$ must lie in the Jacobson radical of $\JJ$.  Hence $1+n$ is invertible, implying that the triangle
    \[
    A\otimes P\to B\otimes P\to C\otimes P\to\Sigma(A\otimes P)
    \]
    is split in $\ul{\Cs}$.  Thus $0\to \varphi(A\otimes P)\to\varphi(B\otimes P)\to\varphi(C\otimes P)\to 0$ is exact, and since $\varphi(P)\otimes(-)$ is faithful and exact, this forces $0\to\varphi(A)\to\varphi(B)\to\varphi(C)\to0$ to also be exact, as desired. 
\end{proof}

\section{Saturation conjecture of Coulembier-Flake}

For this section, $\Bbbk$ denotes an algebraically closed field of characteristic $p>0$.  In this section, we establish a generalisation of \cite[Conj.~5.4.2(1)]{CF} in order to construct semi-exact functors in positive characteristic.

Let $G=SL_2$, and write 
\[
U=\left\{\begin{pmatrix}1 & *\\0 & 1\end{pmatrix}\right\}\leq G
\]
for a root subgroup with normalising Borel subgroup $B$. 
 For $\underline{\lambda}=(\lambda_1,\dots,\lambda_r)\in\A^r(\Bbbk)$, we obtain a morphism $\psi_{\underline{\lambda}}:C_p^r\to U$, where $C_p$ denotes the cyclic group of order $p$, given by:
\[
\psi_{\underline{\lambda}}(a_1,\dots,a_r)=\begin{pmatrix}
    1 & \sum\limits_{i=1}^ra_i\lambda_i \\ 0 & 1
\end{pmatrix}.
\]
Let $\G_{a(s)}$ denote the $s$th Frobenius kernel of $\G_a$.  For $\underline{\mu}=(\mu_0,\dots,\mu_{s-1})\in\A^s(\Bbbk)$, we have a morphism
\[
\sigma_{\underline{\mu}}:\G_{a(s)}\to U, \ \ \ t\mapsto \sum\limits_{i=0}^{s-1}\mu_it^{p^i}.
\]
Write $U_{\underline{\mu}}=\operatorname{im}(\sigma_{\underline{\mu}})$ for the image subgroup scheme.  Then $\sigma_{\underline{\mu}}$ is an embedding if and only if $\mu_0\neq0$.

Putting them together, for a choice of $(\ul{\lambda},\ul{\mu})$, we obtain a map
\[
\psi_{\ul{\lambda}}\times\sigma_{\ul{\mu}}:C_p^r\times\G_{a(s)}\to U,
\]
and this map is an embedding if and only if the $\lambda_i$ are linearly independent over $\F_p$ and $\mu_0$, and the image depends only on $[\ul{\lambda}:\ul{\mu}]\in\mathbb{P}^{r+s-1}$.  

For $\ul{\lambda},\ul{\mu}$ which define an embedding, we set
\[
H_{\underline{\lambda},\underline{\mu}}=E_{\underline{\lambda}}\times U_{\underline{\mu}}\cong C_p^r\times \mathbb{G}_{a(s)},
\]
which is a subgroup of $U$.

Let $U_n$ denote the $n$th Frobenius kernel of $U$.  We have natural restriction functors
\[
R_U:\Rep G\to \Rep U, \ \ R_{\underline{\lambda},\underline{\mu}}:\Rep G\to \Rep H_{\underline{\lambda},\underline{\mu}}, \ \ R_{U_n}:\Rep G\to \Rep U_n.
\]

Write $\operatorname{Tilt}$ for the category of tilting modules for $SL_2$.  This is a pseudo-tensor category, and the tensor ideals are exactly $I_0=\operatorname{Tilt}\supset I_1\supset I_2\supset\cdots$, where $I_n=\langle\operatorname{St}_n\rangle$ (see \cite{C}).  Here $\operatorname{St}_n$ denotes the $n$th Steinberg module of $SL_2$.  The main conjecture of \cite{CF} can be stated as follows.

\begin{conj}\label{conj}
    Let $V\in\Rep E_{\underline{\lambda}}$.  Then $V\otimes R_{\ul{\lambda},0}(\operatorname{St}_{r-1})\in R_{\ul{\lambda},0}(I_{r-1})$.
\end{conj}

We aim to prove Conjecture \ref{conj} by proving the following stronger theorem.

\begin{thm}\label{thm CF saturation}
    Let $V\in\Rep H_{\underline{\lambda},\underline{\mu}}$.  Then $V\otimes R_{\underline{\lambda},\underline{\mu}}(\operatorname{St}_{r+s-1})\in R_{\underline{\lambda},\underline{\mu}}(I_{r+s-1})$.
\end{thm}

In the case $s=0$, we obtain the following as a corollary:

\begin{cor}\label{cor 3.3}
      Let $V\in\Rep (U_s)$.  Then $V\otimes R_{U_s}(\operatorname{St}_{s-1})\in R_{U_n}(I_{s-1})$.
\end{cor}

We now turn to the proof of Theorem \ref{thm CF saturation}.  Write $\operatorname{Dist}(U)=\Bbbk[u_0,u_1,\dots,]/(u_0^p,u_1^p,\dots)$ for the distribution algebra of $U$, where $u_i$ denotes the $(p^i)$-th divided power of a generator $u$ of $\operatorname{Lie}(U)$.  Recall that $\Rep(U)$ is equivalent to the category of $\operatorname{Dist}(U)$-modules $V$ such that for $u_iV=0$ for $i\gg 0$ (\cite{F}).  Thus $\Rep(U_n)$ is the category of modules over the subalgebra $\Gamma_n=\Bbbk[u_0,\dots,u_{n-1}]/(u_0^p,\dots,u_{n-1}^p)$. We have a natural, exact (although not monoidal) functor $\alpha:\Rep(U_n)\to\Rep(U)$ such that $\alpha\circ R_u=\operatorname{id}_{\Rep(U_n)}$, and $u_i$ acts by zero on $\alpha(V)$ for $i\geq n$.
 
\subsection{First step: formulating a sufficient condition}     Let $n=r+s$.  We have a commutative diagram 
\[
\xymatrix{\Rep(U_n) \ar[r]^{\alpha}\ar[dr]^{\beta} & \Rep(U)\ar[d]^\rho\\
& \Rep H_{\underline{\lambda},\underline{\mu}}}
\]
Here $\rho$ is the natural restriction functor.  

\begin{lemma}
    The functor $\beta=\rho\circ\alpha$ is an equivalence of abelian categories.
\end{lemma}
\begin{proof}
    Indeed, first observe that $\operatorname{St}_n$ restricts to the free module of rank one over both $U_n$ and $H_{\underline{\lambda},\underline{\mu}}$ (this can be shown using the polynomial model of the Steinberg representation, see e.g.~\cite[\S 5.3]{CF} and \cite[Lem.~3.2]{CF2}).  Moreover, $\alpha(R_{U_n}(\operatorname{St}_n))=R_U(\operatorname{St}_n)$ by the explicit polynomial realisation of $\operatorname{St}_{n}$.  Thus $\beta$ takes a free module of rank one to a free module of rank one.  Since $\beta$ is exact and faithfully, and fully faithful on projectives, it follows it is a Morita equivalence.
\end{proof}

Now, let $V_{r,s}\in\Rep H_{\underline{\lambda},\underline{\mu}}$, and let $V_{U_n}\in\Rep U_n$ be such that $\beta(V_{U_n})\cong V_{r,s}$.  Finally, let $V=V_U:=\alpha(V_{U_n})$, so that $\rho(V)=V_{r,s}$ and $V|_{U_n}=V_{U_n}$. Our goal is to prove the following.
\begin{thm}\label{thm stronger}
    We have $V\otimes R_U(\operatorname{St}_{n-1})\in R( I_{n-1})$.
\end{thm}

It is easy to see that Theorem \ref{thm stronger} implies both Theorem \ref{thm CF saturation}, and thus also Conjecture \ref{conj} and Corollary \ref{cor 3.3}.

\subsection{Second step: studying the fibre of an induced module}

Consider $M=\Ind_{U}^{G}V.$  This is naturally a $G$-equivariant module over $\Ind_U^G(\Bbbk)=\Bbbk[\A^2\setminus\{(0,0)\}]=\Bbbk[x,y]$.  

\begin{lemma}
    $M$ is a free $\Bbbk[x,y]$-module of rank $\dim V$.
\end{lemma}

\begin{proof}
    Indeed, because $\operatorname{Dist}(U_n)$ is self-injective, there exists an exact sequence of $\operatorname{Dist}(U_n)$-modules $0\to V\to F_1\to F_2$, where $F_1$ and $F_2$ are free.  However $\operatorname{Dist}(U_n)\cong R_{U_n}(\operatorname{St}_{n})$ as a module over itself, so $F_i\cong\operatorname{R}_{U_n}(\operatorname{St}_{n})^{\oplus a_i}$.  Applying left exactness of induction and the isomorphisms
    \[
    \Ind_U^G(\operatorname{St}_n)^{\oplus a_i}\cong \operatorname{St}_n^{\oplus a_i}\otimes\Bbbk[x,y],
    \]
    yields short exact sequences of $\Bbbk[x,y]$-modules
    \[
    0\to M\to \operatorname{St}_n^{\oplus a_1}\otimes\Bbbk[x,y]\to \operatorname{St}_n^{\oplus a_2}\otimes\Bbbk[x,y].
    \]
    Because $\Bbbk[x,y]$ has global dimension 2, this implies $M$ is projective by a standard syzygy argument.
\end{proof}

Write $\m=(x,y)$ for the maximal ideal of $(0,0)\in\A^2$, which is $G$-invariant. Then 
\[
W:=M/\m M
\]
will be a $G$-module of dimension $\dim V$.  

For the following, recall that $u_i$ are the divided powers of a generator $u$ of $\operatorname{Lie}(U)$.
\begin{lemma}
      We have $u_i(W)=0$ for $i\geq n$.
\end{lemma}

\begin{proof}
    Indeed, let $\A^1=\Spec \Bbbk[x]\sub\A^2$ be the line fixed pointwise by $U$, and let \linebreak $\ol{M}=M\otimes_{\Bbbk[x,y]}\Bbbk[x]$.  Then $M$ is a $B$-equivariant vector bundle on $\A^1$ with fibre $W$ at 0 and fibre $V$ at 1.  For a point $a\in\A^1\setminus \{0\}$, we have 
    \[
    \begin{pmatrix}
        a & 0\\0 & a^{-1}
    \end{pmatrix}\cdot 1=a
    \]
    Thus the action of $U$ on the fibre of $\ol{M}$ at $a$ is obtained by conjugating the action at $1$ by $\begin{pmatrix}
        a & 0\\0 & a^{-1}
    \end{pmatrix}$.  The operators $u_i$ are eigenvectors for this conjugation.  Since $u_i$ acts by 0 on the fibre at 1 for $i\geq n$, these operators must act by zero on the fibre of $\ol{M}$ at every nonzero $a$.  
    
    But $U$ acts trivially on the underlying space $\A^1$, so the operators $u_i$ will act $\OO_{\A^1}$-linearly on $\ol{M}$.  Since $\ol{M}$ is a vector bundle and these operators act trivially away from 0, they must act trivially everywhere, including the fibre $W$ at 0, completing the proof.
\end{proof}

In the following, write $L_a$ for the irreducible $SL_2$ module of highest weight $a\in\N$.

\begin{cor}\label{cor comp fact X}
    If $L_a$ is a composition factor of $W$, then $a<p^{n}$.
\end{cor}

\begin{proof}
    Indeed, if $a\geq p^n$ then $u_{n+1}$ must act nontrivially by the Steinberg tensor product theorem.
\end{proof}

Corollary \ref{cor comp fact X} and \cite[Lem.~3.3]{BEO} now imply:
\begin{cor}\label{cor}
    $W\otimes\operatorname{St}_{n-1}$ is a tilting module.
\end{cor}

Now we consider $Z:=M\otimes\operatorname{St}_{n-1}$.  It is easy to see that $T:=Z/\m Z\cong W\otimes\operatorname{St}_{n-1}$, which is tilting by Corollary \ref{cor}. 

\subsection{Third step: producing a $G$-equivariant splitting of $Z\to T$}

Observe that for all $i$, we have a natural morphism of $G$-modules
\[
\m^i/\m^{i+1}\otimes T\to \m^iZ/\m^{i+1} Z.
\]
This must be an isomorphism since $M\otimes \operatorname{St}_{n-1}$ is a free $\Bbbk[x,y]$-module.  On the other hand, we have $\m^i/\m^{i+1}\cong \nabla(i)$, where $\nabla(i)$ denotes the costandard module of highest weight $i$.  Thus $\m^i/\m^{i+1}\otimes T$ has a filtration by costandard modules. 

Now consider the short exact sequences:
\[
0\to \m^i Z/\m^{i+1}Z\to Z/\m^{i+1} Z\to Z/\m^i Z\to 0
\]
The object $\m^i Z/\m^{i+1}Z\cong\m^i/\m^{i+1}\otimes T$ has a filtration by costandard modules, and thus
\[
\Ext^1(T,\m^iZ/\m^{i+1}Z)=0 \text{ for all }i,
\]
since $T$ is tilting.  For $i=1$ this provides a splitting of $T$ off $Z/\m^2Z$, and then at $i=2$ a splitting off $Z/\m^3Z,$ and so on.

Working inductively, for each $i$ we obtain a splitting of $T$ off $Z/\m^{i+1}Z$, which is compatible with the splitting off $Z/\m^{i}Z$.  Therefore we obtain a map of vector spaces 
\[
T\to\lim\limits_{\leftarrow}Z/\m^iZ.
\]
Moreover, this map is equivariant for the natural action of $G(\Bbbk)$. 

Recall we have a short exact sequence of $G$-equivariant $\Bbbk[x,y]$-modules
\[
0\to Z\to\operatorname{St}_{n-1}\otimes \operatorname{St}_n^{\oplus a_1}\otimes\Bbbk[x,y]\to \operatorname{St}_{n-1}\otimes \operatorname{St}_n^{\oplus a_1}\otimes\Bbbk[x,y].
\]
If we complete this sequence with respect to $\m$ it will remain exact (\cite[Lem.~10.97.1]{Stacks}), and thus we have embeddings of $G$-modules
\[
\xymatrix{
 & Z \ar@{^{(}->}[d] \ar@{^{(}->}[r] & \operatorname{St}_{n-1}\otimes \operatorname{St}_n^{\oplus a_1}\otimes\Bbbk[x,y] \ar@{^{(}->}[d] \\ T\ar@{^{(}->}[r] & \hat{Z}\ar@{^{(}->}[r] & \prod\limits_{a=0}^\infty \operatorname{St}_{n-1}\otimes \operatorname{St}_n^{\oplus a_1}\otimes\nabla(a).
}
\]
This diagram is moreover a pullback.  We now observe that $T$ has a nonzero map to $\operatorname{St}_{n-1}\otimes \operatorname{St}_n^{\oplus a_1}\otimes\nabla(a)$ for only finitely many $a$, since the socle of $\nabla(a)$ is $L(a)$.  Hence $T$ lands in $\operatorname{St}_{n-1}\otimes \operatorname{St}_n^{\oplus a_1}\otimes\Bbbk[x,y]$ under the composition of the bottom horizontal embeddings, and hence must land in $Z$.  We may conclude:

\begin{lemma}
    The map $Z\to Z/\m Z=T$ splits as a $G$-module.
\end{lemma}

\subsection{Fourth step: obtaining a global $G$-equivariant splitting}

Consider the $G$-equivariant map induced by the $\Bbbk[x,y]$-module action:
\[
\varphi:\Bbbk[x,y]\otimes T\to Z.
\]
We claim this is an isomorphism.  Indeed, it is an isomorphism in a neighborhood of 0 since $T$ maps isomorphically onto the fibre there.  Since the complement of 0 is a homogeneous $G$-space, and $\varphi$ is $G$-equivariant, it must be a global isomorphism.  

Taking the fibre at the base point $eU\in G/U(\Bbbk)$, we obtain an isomorphism of $U$-modules $V\otimes R(\operatorname{St}_{n-1})\cong T$, which finishes the proof of Theorem \ref{thm 3}.

Recalling that $T=W\otimes \operatorname{St}_{n-1}$, we additionally obtain the following corollary.

\begin{cor}\label{cor 2}
    Let $V_{r,s}\in\Rep(H_{\underline{\lambda},\underline{\mu}})$.  Then there exists $W\in\Rep(SL_2)$, all of whose weights are $p^n$-restricted, such that 
    \[
    V_{r,s}\otimes R_{\underline{\lambda},\underline{\mu}}(\operatorname{St}_{n-1})\cong R_{\underline{\lambda},\underline{\mu}}(W\otimes\operatorname{St}_{n-1}).
    \]
\end{cor}

\section{Proof of Theorem \ref{thm intro}}

\subsection{The Coulembier-Flake functor}

In \cite{CF}, an LSM functor $\Rep C_p^r\to\Ver_{p^r}$ is constructed (subject to the validity of Conjecture \ref{conj}), which depends on the choice of a sufficiently generic point $[\underline{\lambda}]\in\mathbb{P}^{r-1}$.  We now generalise this construction to a quasi-elementary abelian group $E=C_p^r\times\G_{a(s)}$, to produce a functor 
\[
\varphi_{\underline{\lambda},\underline{\mu}}:\Rep E\to\Ver_{p^n},
\]
where $n=r+s$, which will depend on the choice of $[\underline{\lambda}:\ul{\mu}]\in\mathbb{P}^{r+s-1}$.  The construction we give follows the ideas from Section \ref{sec abel env}.

Let $S:=R_{\underline{\lambda},\underline{\mu}}(\operatorname{St}_{n-1})$, and let $\II$ denote the additive tensor ideal in $\ul{\Rep E}$ generated by $S$.  Then by Theorem \ref{thm CF saturation}, $\II$ realises a point $\PP\in\operatorname{cl}\Sp\ul{\Rep E}$ which is described in \cite[Cor.~4.2.7]{CF3}.  Thus we may apply Theorem \ref{thm abelian envelope semi exact} to obtain an abelian envelope
\[
\ul{\Rep E}/\PP^{max}\to \Ds.
\]
Thus we obtain a functor
\[
\operatorname{Tilt}SL_2\to\Rep E\to \ul{\Rep E}/\PP^{max}\to \Ds,
\]
and we know that $\operatorname{St}_{n-1}$ is not killed while $\operatorname{St}_{n}$ goes to 0.  Thus, writing $I_n$ for the tensor ideal of $\operatorname{Tilt}SL_2$ generated by $\operatorname{St}_n$, and using the classification of tensor ideals in $\operatorname{Tilt}SL_2$ \cite{C}, we obtain a faithful LSM functor
\[
\operatorname{Tilt}SL_2/I_n\to \Ds,
\]
which therefore factors through a tensor functor from the abelian envelope of $\operatorname{Tilt}SL_2/I_n$, which is $\Ver_{p^n}$:
\[
F:\Ver_{p^n}\to\Ds.
\]

\begin{lemma}
    The functor $F$ is an equivalence of categories.
\end{lemma}
\begin{proof}
    Because $\Ver_{p^n}$ is incompressible \cite{BEO}, it suffices to prove that $F$ is a surjective tensor functor.  By the theory of abelian envelopes explained in Section \ref{sec abel env}, the indecomposables of $\II$ will make up all of the projective indecomposables of $\Ds$.  However $\II$ is in the image of the functor $\operatorname{Tilt}SL_2\to\ul{\Rep E}$, and so the image of $F$ must contain all projectives in $\Ds$.  From here it is straightforward to show that $F$ is surjective. 
\end{proof}

We have thus constructed the desired LSM functor
\[
\varphi_{\underline{\lambda},\underline{\mu}}:\Rep E\to\Ver_{p^n}.
\]
  
\begin{thm}\label{thm 3}
   The functor $\varphi_{\underline{\lambda},\underline{\mu}}$ is semi-exact.
\end{thm}

\begin{proof}
    Indeed, this follows from Theorem \ref{thm abelian envelope semi exact}.
\end{proof}

\subsection{Finite group schemes and quasi-elementary subgroups}

Let $G$ be a finite group scheme, and let $\Cs=\Rep G$.  Recall that a quasielementary group scheme $E$ is one isomorphic to $\C_p^r\times\G_{a(s)}$ for some $r,s\in\N$. For every point $x\in\operatorname{cl}\Sp\ul{\Cs}$, there exists a minimal quasi-elementary subgroup $E\leq G$ such that $x$ lies in the image of $\Sp\ul{\Rep E}\to\Sp\ul{\Cs}$ \cite{FP}.

\subsection{Finite supergroup schemes in characteristic 0}\label{section super}

If $\Bbbk$ is characteristic 0 and $\Cs$ is a finite symmetric tensor category, then by \cite{Del}, $\Cs\simeq\Rep(H\ltimes\Pi V,\varepsilon)$ for a finite group $H$, an $H$-module $V$, and $\varepsilon\in Z(H)$ such that $\varepsilon^2=1$, and $\varepsilon$ acts by $-1$ on $V$.  We write $G:=H\ltimes\Pi V$.  In this case, the works of \cite{BKN}, \cite{HV}, and \cite{DK} give the Balmer spectrum as
\[
\Sp\ul{\Cs}=\mathbb{P}(V)/H.
\]
Thus to a point $x\in\operatorname{cl}\Sp\ul{\Cs}$ we will associate an element $x\in V$, which is well-defined up to non-zero scalar and $H$-orbit.

\subsection{The proof of Theorem \ref{thm intro}}

We now give a proof of Theorem \ref{thm intro}.  Therefore we either assume that $\Bbbk$ is characteristic 0 and $\Cs=\Rep(G,\varepsilon)$ with $G=H\ltimes\Pi V$ as in the setup of Section \ref{section super}, or that $\Bbbk$ is characteristic $p>0$ and $\Cs=\Rep G$ for a finite group scheme $G$.

\begin{proof}[Proof of Theorem \ref{thm intro}]
    Let $\varphi:\Cs\to\Ds$ be a surjective semi-exact functor.  If $\varphi$ is exact, then $\Ds$ must be finite by \cite[Prop.~1.8.19]{EGNO}, and thus we now assume that $\varphi$ is not exact.  Then by Theorem \ref{thm m_phi maximal}, $\varphi$ determines a point $\KK_{\varphi}\in\operatorname{cl}\Sp\ul{\Cs}$.  
    
    If $\Bbbk$ is characteristic 0, then by Section \ref{section super}, $\KK_{\varphi}$ lifts to a non-zero element $x\in V$, and in this case we let $E=\Pi\langle x\rangle$ denote the purely odd, one-dimensional supergroup generated by $x$. Then we set $S:=\Ind_{E}^{G}\ol{S},$ where $\ol{S}$ denotes the trivial $E$-module.
    
    If $\Bbbk$ is positive characteristic, then $\KK_{\varphi}$ arises from a quasi-elementary abelian subgroup $E$ of rank $n$ such that it restricts to a generic point $([\ul{\lambda}],[\ul{\mu}])$ in $\Sp\ul{\Rep E}$.  In this case, let $\ol{S}=R_{\ul{\lambda},\ul{\mu}}(St_{n-1})\in\Rep E$, and set $S:=\Ind_{E}^{G}\ol{S}$.
    \begin{lemma}\label{lemma supp of S}
        The support of $S$ is exactly the point $\KK_{\varphi}$.
    \end{lemma}

    \begin{proof}
        In characteristic $p>0$ we may apply \cite[Prop.~9.5(2)]{BIKP}.  In characteristic 0, the induced module is given by $\bigoplus\limits_{g\in G}\Lambda(V/\langle gx\rangle)$, which has support precisely the $G$-orbit of $x$, as desired.
    \end{proof}

    Recall that for $M\in\Cs$, we have the identity 
    \[
    M\otimes\Ind_{E}^{G}\ol{S}\cong \Ind_{E}^{G}(\Res_E(M)\otimes \ol{S}).
    \]
    In all cases we have that $\ol{S}$ generates an additive tensor ideal in $\ul{\Rep E}$ with finitely many indecomposables, and thus $S$ must do the same in $\ul{\Cs}$.  Write $\II$ for the additive tensor ideal generated by $S$.  Then since $S$ is supported exactly at $\m_{\varphi}$ by Lemma \ref{lemma supp of S}, $\II$ realises $\KK_{\varphi}$ in the sense of Section \ref{sec abel env}. By Theorem \ref{thm abelian envelope semi exact} there exists an abelian envelope
    \[
    \ul{\Cs}/\KK_\varphi^{max}\to \Ds'.
    \]

    \begin{lemma}\label{lemma varphi kills k_x}
        $\ker(\varphi)=\KK_\varphi^{max}$.
    \end{lemma}

    \begin{proof}
        Because $\varphi(M)=0$ if and only if $\KK_{\varphi}$ does not lie in the support of $M$, it is clear that $\ker(\varphi)\sub\KK_\varphi^{max}$.  We must prove the opposite inclusion.

        Let $f:X\to Y$ be a morphism in $\KK_\varphi^{max}$.  Then $f\otimes S:X\otimes S\to Y\otimes S$ must lie in the categorical radical $\operatorname{Rad}(\II)$ of the finite additive category $\II$.  Since $\II$ has finitely many indecomposables, there exists $N>0$ such that $\operatorname{Rad}(\II)^N=0$.  We claim that $f^{\otimes N}\otimes S=0$.  Indeed, it may be factored as
        \[
        X^{\otimes N}\otimes S\xto{f\otimes S}Y\otimes X^{\otimes (N-1)}\otimes S\to\cdots Y^{\otimes N}\otimes S,
        \]
        and thus by assumption is 0.  It follows that $S\otimes F(f)^{\otimes N}=0$ in $\Ds$.  But since $\Ds$ is a tensor category, this forces $F(f)=0$, completing the proof.
    \end{proof}

    It now follows that $\varphi$ defines a faithful symmetric monoidal functor $:\ul{\Cs}/\KK_\varphi^{max}\to\Ds$, and thus by the universal property we obtain a surjective, exact tensor functor $F:\DD'\to\DD$ fitting into the following diagram:
\[
\xymatrix{
\Cs\ar[dr]^{\varphi}\ar[r] \ar@/^2pc/[rr]^{\tilde{\varphi}} & \ul{\Cs}/\KK_\varphi^{max}\ar[r]\ar[d] & \Ds'\ar[dl]_{F} \ar[d]^{\omega'} \\
& \Ds\ar[r]_{\omega} & \Ver_{p^\infty}}
\]
   Since $F$ is surjective and exact, and $\Ds'$ is finite, $\Ds$ must also be finite.  Since $\Ds,\Ds'$ are finite, there exists exact tensor functors $\omega,\omega'$ as above, by the main theorem of \cite{CF2}.  By \cite[Thm.~4.7]{CF2}, $\omega\circ F$ must be isomorphic to $\omega'$, and thus 
   \[
   \omega'\circ\tilde{\varphi}\simeq \omega\circ\varphi.
   \]
    It follows that if $\Ds'=\Ver_{p^\infty}$, then $\varphi\simeq\tilde{\varphi}$, which proves the injectivity of the map in Conjecture \ref{conj 1}.  Surjectivity follows from the existence of $\tilde{\varphi}$.
\end{proof}

    \bibliographystyle{amsalpha}

\textsc{\footnotesize A.S.: School of Mathematics and Statistics, UNSW Sydney, NSW 2052,
		Australia}

	\textit{\footnotesize Email address:} \texttt{\footnotesize xandersherm@gmail.com}
\end{document}